\documentclass[12pt]{amsart}
\usepackage{amsmath,amsthm,amssymb}
\usepackage{lscape}
\usepackage[all]{xy}
\usepackage{pifont}
\usepackage{bbding}
\usepackage{url}

\newtheorem{thm}{Theorem}[section]
\newtheorem{lem}[thm]{Lemma}
\newtheorem{dfn}[thm]{Definition}

\newtheorem{rem}[thm]{Remark}

\newtheorem{hyp}[thm]{Hypothesis}
\newenvironment{prf}{{\bf Proof.}}{\hfill $\diamond$}

\def\dim{ {\rm dim } }

\def\Aut{ {\rm Aut } }

\def\Sp{ \text{\rm Sp} }

\def\PSp{ \text{\rm PSp} }
\def\SL{ \text{\rm SL} }
\def\PSL{ \text{\rm PSL} }
\def\PGL{ \text{\rm PGL} }

\def\SU{ \text{\rm SU} }
\def\PSU{ \text{\rm PSU} }

\begin{document}

\title{On maximal embeddings of finite quasisimple groups}

\author{Gerhard Hiss}

\address{Lehrstuhl f{\"u}r Algebra und Zahlentheorie, RWTH Aach\-en University,
52056 Aach\-en, Germany}

\email{gerhard.hiss@math.rwth-aachen.de}



\dedicatory{Dedicated to the Memory of my Friend and Colleague Kay Magaard}

\date{\today}

\begin{abstract}
If a finite quasisimple group~$G$ with simple quotient~$S$ is embedded into a 
suitable classical group~$X$ through the smallest degree of a projective 
representation of~$S$, then $N_X( G )$ is a maximal subgroup of~$X$, up to 
two series of exceptions where~$S$ is a Ree group, and four exceptions 
where~$S$ is sporadic. 
\end{abstract}

\maketitle

\section{Introduction} \label{Introduction}

With his seminal paper~\cite{Asch}, Michael Aschbacher initiated the programme
of classifying the maximal subgroups of the finite classical groups. He 
introduced eight classes of geometric subgroups of such a classical group
and showed that a maximal subgroup lies in one of these classes or is what is 
now called a subgroup of type~$\mathcal{S}$. After the monumental works of 
Kleidman and Liebeck \cite{KL} and of Bray, Holt and Roney-Dougal \cite{BHRD}, 
it remains to determine the maximal subgroups of type~$\mathcal{S}$ for 
classical groups of degree at least~$13$.

If $X$ is a quasisimple classical group, a maximal subgroup of type 
$\mathcal{S}$ of~$X$ is of the form $N_X( G )$, where~$G$ is a quasisimple 
group, acting absolutely irreducibly on the vector space underlying~$X$ (and 
satisfies further conditions). One may now ask conversely: Which absolutely 
irreducible faithful representations of a quasisimple group~$G$ into a classical 
group $X$ give rise to maximal subgroups of~$X$ of type~$\mathcal{S}$? 
By \cite[Theorem~$7.5$(b), Example~$7.6$]{HiHuMa}, most such representations
(in a precise quantitative sense) are imprimitive and thus, in general, do not
yield maximal subgroups of~$X$, as their image is contained in the
stabilizer of the imprimitivity decomposition of the underlying vector space.
In \cite[Section~$3$]{Kay}, Kay Magaard has answered the above question for the 
simple Mathieu group $G = M_{11}$. Surprisingly enough, only the two 
$5$-dimensional representations of~$M_{11}$ over the field with three elements 
give maximal subgroups of~$X$. While this somewhat unexpected behaviour may be 
due to the small size of the group, another feature of this example is not. The 
degree~$5$ representations of $M_{11}$ are its representations of smallest 
possible degree greater than~$1$. In this note we show that this is a general 
phenomenon. 
Let~$G$ be a quasisimple group and let~$n$ denote the smallest degree of 
a non-trivial projective representation of the simple group $G/Z(G)$. Then, up
to two series of exceptions and four single cases the following holds: 
If $X$ is a quasisimple classical 
group of degree~$n$ properly containing~$G$ and minimal with this property, then
$N_X( G )$ is maximal in~$X$.

\section{Minimal projective degrees} \label{MinProjDeg}

Unless otherwise stated, the groups in our paper will be finite. For a prime 
power~$q$ we write $\mathbb{F}_q$ for the finite field with $q$ elements and 
$\overline{\mathbb{F}}_q$ for its algebraic closure. If $H$ and~$G$ are groups,
we write $H \preceq G$, if $H$ is isomorphic to a subgroup of~$G$, and we 
write~$G'$ for the derived subgroup and $G^\infty$ for the last term in the 
derived series of~$G$. We will 
make use of the following well-known result about quasisimple groups, which 
follows from the Three Subgroups Lemma.

\begin{lem}
\label{QSExtensions}
Let $G$ be a perfect group, $Z \leq Z(G)$ such that $G/Z$ is quasisimple. 
Then~$G$ is quasisimple.
\end{lem}

\smallskip

\noindent
The following definitions are taken from \cite[(5.3.1),(5.3.2)]{KL}. 

\begin{dfn}[Kleidman-Liebeck, \cite{KL}]
\label{RpG}
{\rm
Let $G$ be a finite group.

(a) If $p$ is a prime, set
$$R_p(G) := \min\{ 0 \neq n \in \mathbb{N} \mid G \preceq \PGL_n( \overline{\mathbb{F}}_p ) \}.$$
Notice that
$$R_p(G) = \min\{ 0 \neq n \in \mathbb{N} \mid G \preceq \PGL_n( F ), 
F \text{\ a field of characteristic\ } p \}$$
(see \cite[(5.3.2)]{KL}).

(b) We also set
$$R(G) := \min\{ R_p( G ) \mid p \text{\ a prime}\}.$$
}
\end{dfn}

\smallskip

\noindent 
We collect a few properties of these invariants needed later on.

\begin{rem} \label{PropertiesR}
{\rm 
Let $S$ be a non-abelian simple group with universal covering group $\hat{S}$ 
and let $p$ be a prime.

(a) Let $F$ be a field of characteristic $p$ and $\rho$ an $F$-representation of 
$\hat{S}$ of degree $R_p(S)$ with non-trivial image. Then $\rho$ is absolutely 
irreducible. 

(b) If $U \lneq S$ is a proper subgroup, then $R(S) < [S\colon\!U]$.

(c) Let $G$ and $K$ be quasisimple with $G \leq K$. Then $R_p( K/Z( K ) ) \geq 
R_p( G/Z( G ) )$. In particular, $R( K/Z( K ) ) \geq R( G/Z( G ) )$.
}
\end{rem}

\smallskip

\noindent
If $S$ is a non-abelian simple group, then $R_p( S )$ is the smallest degree of
a non-trivial projective representation of~$S$ over a field of characteristic~$p$
(see \cite[Proposition~$5.3.1$(ii)]{KL}). Using this, part~(a) of 
Remark~\ref{PropertiesR} is easily verified. Using in addition that the universal
covering group of $K/Z( K )$ contains a covering group of $G/Z( G )$ by
Lemma~\ref{QSExtensions}, we obtain part~(c). To prove~(b), consider the 
permutation representation of~$S$ on the cosets of~$U$.

If $S$ is a sporadic simple group, $R(S)$ can be determined from~\cite{Jansen}.
If $S = A_m$ is an alternating group of degree $m \geq 5$, then $R( A_m ) < 5$
for $5 \leq m \leq 8$, and $R( A_m ) = m - 2$ for $m \geq 9$ 
(see \cite[Proposition~$5.3.7$]{KL}). The value $R(S)$ for a simple group~$S$ 
of Lie type can be determined from 
\cite[Propositions~$5.4.13$,~$5.4.15$, Corollary~$5.4.14$]{KL} and~\cite{MAtl}.

\section{The main result}

Keep the notation of the previous section. We shall consider the following 
hypothesis.

\begin{hyp}
\label{BasicHypothesis}
{\rm
Let $S$ be a finite non-abelian simple group and put $n := R( S )$. Assume that
$n \geq 5$. Let $p$ be a prime such that $R( S ) = R_p( S )$. Let $G$ be a 
covering group of~$S$ such that~$G$ has a faithful representation 
$\rho$ of degree~$n$ over $\overline{\mathbb{F}}_p$. Let $q$ be a power of~$p$ 
minimal with the property that~$\rho$ is realizable over~$\mathbb{F}_q$.

Let $V$ be an $\mathbb{F}_qG$-module affording~$\rho$ and identify~$G$ with 
$\rho( G ) \leq \SL( V )$. Finally, let $X \leq \SL( V )$ denote the smallest 
quasisimple classical group containing~$G$.
}
\end{hyp}

Let us clarify what we mean by the smallest quasisimple classical group 
containing~$G$. If~$G$ stabilizes a non-degenerate quadratic form on~$V$, 
then~$X$ is the commutator subgroup of the full isometry group of this form in 
$\SL( V )$. Suppose that $G$ does not stabilize any non-degenerate quadratic 
form on~$V$. If~$G$ stabilizes a non-degenerate symplectic or hermitian form
on~$V$, then~$X$ is the full isometry group of this form in $\SL(V)$. Otherwise 
$X = \SL( V )$. Notice that $X$ is uniquely determined by~$G$ (see
\cite[Lemma~$1.8.8$]{BHRD}). Moreover,~$X$ is quasisimple, as $n \geq 5$. 
Finally, $X$ is isomorphic to one of $\Omega_n( q )$ (where $n$ is 
odd), $\Omega^\pm_n( q )$ (where $n$ is even), $\Sp_n( q )$ (where $n$ is even), 
$\SU_n( q_0 )$ (where $q = q_0^2$) or $\SL_n( q )$.

\begin{lem}\label{UpperLowerGroup}
Assume {\rm Hypothesis~\ref{BasicHypothesis}}. Let $K$ be quasisimple such that 
$G \leq K \leq X$, and put $T := K/Z( K )$. Then $R_p( T ) = R_p( S )$ and 
$R( T ) = R( S )$. If~$T$ is a group of Lie type of characteristic~$\ell$, then 
$\ell = p$ unless $T = \Omega^+_8( 2 )$ or $T = {^2\!F}_4(2)'$. 

If $T$ is a classical group, then $K = X$, unless $T = \Omega^+_8( 2 )$ and~$p$
is odd. In the latter case, $G = K$. 

If $S$ and $T$ are exceptional groups of Lie type, or if $T = {^2\!F}_4(2)'$,
then $G = K$, unless $(G,K) = ({^2G}_2( q ), G_2( q ))$, $q = 3^{2m+1}$
with $m \geq 1$, or $(G,K) = ({^2\!F}_4( q )', F_4( q ))$, $q = 2^{2m+1}$
with $m \geq 0$.
\end{lem}
\begin{prf}
The embedding of~$K$ into~$X$ is a faithful representation of~$K$ of degree~$n$, 
and thus $R_p( T ) \leq n$ by the comment following Remark~\ref{PropertiesR}. 
As $G \leq K$, we conclude from Remark~\ref{PropertiesR}(c) that $n \geq R_p( T ) 
\geq R_p( S ) = R( S ) = n$. It follows that $R_p( T ) = R_p( S )$ and also that 
$n = R_p( T ) \geq R( T ) \geq R( S ) = n$, i.e.\ $R( T ) = R( S )$.

Suppose now that $T$ is a group of Lie type of defining characteristic~$\ell$. 
As $R( T ) \geq 5$, we conclude from 
\cite[Propositions~$5.4.13$,~$5.4.15$, Corollary~$5.4.14$]{KL} and \cite{MAtl}, 
that $R_\ell( T ) = R( T )$. As $R( T ) = R_p(T) = n$ by the first paragraph of 
this proof, the same references also give $\ell = p$, unless 
$T = \Omega^+_8( 2 )$ or $T = {^2\!F}_4(2)'$, and if $T$ is a classical group of 
natural dimension~$d$, then $d = R( T ) = n$. In the latter case we get $K = X$,
unless $T = \Omega^+_8( 2 )$ and~$p$ is odd. Indeed, $K$ and $X$ are then 
quasisimple classical groups of the same characteristic and degree, they 
stabilize the same form by \cite[Lemma~$1.8.8$]{BHRD}, $K \leq X$ and $K$ is not 
realizable over a smaller field than $\mathbb{F}_q$. Suppose now that $T = 
\Omega_8^+( 2 )$ and~$p$ is odd. Then $n = R( T ) = 8$ by~\cite{MAtl}. Suppose 
that $S \lneq T$ and consider the list of maximal subgroups of~$T$ (see 
\cite[p.~$85$]{Atlas}). This implies that $S \preceq S_1$, where~$S_1$ is one of 
the simple groups $\Sp_6(2)$, $A_8$, $A_9$, $\SU_4( 2 )$ or $A_5$. 
But then $R( S_1 ) < 8$, and hence $R( S ) < 8$ by Remark~\ref{PropertiesR}(c), 
a contradiction. Thus $S = T$ and hence $G = K$.

Suppose now that $S$ and $T$ are exceptional groups of Lie type, but that
$T \neq {^2\!F}_4(2)'$. By \cite[Table~$5.4.C$]{KL} we find that~$R( T )$ 
determines the Dynkin type of~$T$. Thus $S$ and~$T$ are of the same
Dynkin type, and either both are untwisted, both are twisted, or one is 
twisted and the other one is untwisted. In the first two alternatives,
the exceptional groups of Lie type $G \leq K$ are of the same twisted 
type and the same characteristic. In this case, we get $G = K$, as
the minimal field of definition of a representation of~$K$ of degree
$R( K )$ determines $K$
(see \cite[Propositions~$5.4.6$(i), $5.4.17$, $5.4.18$, Remark~$5.4.7$]{KL}).
Now suppose that we are in the third alternative. We begin by considering
the Dynkin type $E_6$, so that $\{ S, T \} = \{ {^2\!E}_6( q_1 ), E_6( q_2 ) \}$
where~$q_1$ and~$q_2$ are powers of~$p$. By \cite[Proposition~$5.4.17$]{KL},
we have $q = q_1^2$ and $X = \SU_{27}( q_1 )$ if $S = {^2\!E}_6( q_1 )$,
and we have $q = q_2$ and $X = \SL_{27}( q )$ if $S = E_6( q_2 )$. Both 
alternatives first lead to $q_2 = q_1^2$ and then to a contradiction. The 
cases of Dynkin type~$F_4$ and $G_2$ yield the pairs $(G,K)$ with $G \lneq K$ 
listed in the last statement of the lemma. (The case $G = {^2G}_2(3)'$
is excluded, as $R({^2G}_2(3)') = 2$ and $R(G_2(3)) = 7$.)

Suppose finally that $T = {^2\!F}_4(2)'$. Then $K = T$ and 
$n = 26$. If $G \lneq K = T$, then~$S$ is a proper section of~$T$. But this 
would imply $R( S ) < 26$ by the list of maximal subgroups of~$T$ (see 
\cite[p.~$74$]{Atlas}), a contradiction. Thus $G = K$ and we are done.
\end{prf}

\smallskip

\noindent
We can now state our main theorem.

\begin{thm}\label{MaxEmb}
Assume {\rm Hypothesis~\ref{BasicHypothesis}}. Then one of the following holds:

{\rm (a)} We have $G = X$.

{\rm (b)} The normalizer $N_X( G )$ is a maximal subgroup in $X$.

{\rm (c)} We have $G = {^2G}_2( q )$ with $q = 3^{2m+1}$, $m \geq 1$ 
and $n = 7$. In this case, $X = \Omega_7( q )$ and $G \lneq G_2( q ) \lneq X$ 
for all $q$.

{\rm (d)} We have $G = {^2\!F}_4( q )'$ with $q = 2^{2m+1}$, $m \geq 0$
and $n = 26$. In this case, $X = \Omega_{26}^+( q )$ and $G \lneq F_4( q ) \lneq X$ 
for all $q$.

{\rm (e)} We have $G = J_2$, $n = 6$ and $q = 4$. In this case $X = \Sp_6( 4 )$
and $G \lneq G_2( 4 ) \lneq X$.

{\rm (f)} We have $G = M_{23}$, $n = 11$ and $q = 2$. In this case,
$X = \SL_{11}( 2 )$ and $G \lneq M_{24} \lneq X$.

{\rm (g)} We have $G = 3.\text{\rm Fi}_{22}$, $n = 27$ and $q = 4$.
In this case, $X = \SU_{27}( 2 )$, and $G \lneq 3.{^2\!E}_6(2) \lneq X$.

{\rm (h)} We have $G = \text{\rm Th}$, $n = 248$ and $q = 3$.
In this case, $X = \Omega^+_{248}( 3 )$, and $G \times 2 \lneq E_8( 3 ) 
\times 2 \lneq X$.
\end{thm}
\begin{prf}
If~$S$ is a classical group, we get $G = X$ and hence conclusion~(a) by 
Lemma~\ref{UpperLowerGroup} (applied with $K = G$), unless 
$S = \Omega^+_8( 2 )$, $n = 8$ and $p$ is odd. Suppose in the following that 
conclusion (a) does not hold. Then $N_X(G)$ is a proper subgroup
of~$X$, as~$G$ and~$X$ are quasisimple. Choose a maximal subgroup $L$ of~$X$ 
with $N_X( G ) \leq L$. We consider the possible Aschbacher classes 
containing~$L$. For the definition of Aschbacher classes we follow 
\cite[Subsections $2.2.1$--$2.2.8$, Definition~$2.1.3$]{BHRD}.

By construction,~$V$ is absolutely irreducible. Clearly,~$V$ is tensor 
indecomposable, as $\dim(V) = R_p(G)$. Also,~$V$ is primitive by 
Remark~\ref{PropertiesR}(b). By definition,~$G$ and hence~$L$ do not lie in
Aschbacher class~$\mathcal{C}_5$. It follows that~$L$ lies in one of the 
Aschbacher classes $\mathcal{C}_3$,~$\mathcal{C}_6$,~$\mathcal{C}_8$ 
or~$\mathcal{S}$. Put $K := L^\infty$. Then $G = N_X( G )^\infty \leq K$.
In particular,~$K$ acts absolutely irreducibly on~$V$. But then~$L$ cannot
lie in Aschbacher class~$\mathcal{C}_3$, as in that case~$L^\infty$ does not
act absolutely irreducibly (see \cite[Lemma~$2.2.7$]{BHRD} and the remark 
following \cite[Definition~$2.1.4$]{BHRD}). Suppose that~$L$ 
lies in Aschbacher class $\mathcal{C}_6$. Then $n = r^m$ for some prime~$r$ and 
some positive integer~$m$. As $n \geq 5$ by our assumption, we have $m \geq 3$ 
if $r = 2$ and $m \geq 2$ if $r = 3$. Moreover, $K = U.C$, where~$U$ is a normal 
$r$-subgroup of~$K$ and~$C$ is a quasisimple classical group of degree $2m$ 
(see \cite[$(7.6.1)$]{KL} or \cite[Table~$2.9$]{BHRD}). It follows that 
$G/(U \cap Z( G ))$ embeds into~$C$, and hence~$G$ has a nontrivial irreducible 
representation over a field of characteristic~$r$ of degree at most $2m$. Thus 
$r^m = n = R( S ) \leq 2m$, contradicting $n \geq 5$. Suppose that~$L$ lies in 
Aschbacher class $\mathcal{C}_8$. Then $K$ is a quasisimple classical group of
characteristic~$p$, and Lemma~\ref{UpperLowerGroup} yields the contradiction 
$L \geq K = X \gneq L$.

We conclude that $L$ lies in Aschbacher class $\mathcal{S}$. In particular,~$K$ 
is quasisimple and $L = N_X( K )$. If $G = K$, we get $N_X( G ) = N_X( K ) = L$, 
hence $N_X( G )$ is a maximal subgroup of~$X$ and conclusion~(b) holds. So let 
us assume that $G \lneq K$ in the following and put $T := K/Z( K )$. As $K \leq 
L \lneq X$, Lemma~\ref{UpperLowerGroup} implies that~$T$ is not a classical 
group. If both of~$S$ or~$T$ are exceptional groups of Lie type, then $(G,K)$ 
is as in (c) or~(d), again by Lemma~\ref{UpperLowerGroup}. The structure of~$X$
in these cases follows from \cite[Table~$5.4.$C]{KL} and the fact that
the Tits group ${^2\!F}_4(2)'$ embeds into $\Omega^+_{26}( 2 )$; see
\cite{MAtl}. We may hence assume 
that not both of~$S$ or~$T$ are exceptional groups of Lie type, and 
$T \neq {^2\!F}_4(2)'$.

Suppose that~$T$ is an exceptional group of Lie type. Then~$p$ is the 
defining characteristic of~$T$ by 
Lemma~\ref{UpperLowerGroup}. By the considerations in the previous paragraph and
at the beginning of our proof, we 
may then assume that~$S$ is a sporadic group, an alternating group, or $S = 
\Omega^+_8( 2 )$, $n = 8$ and~$p$ is odd. 
Since $R( T ) \leq 248$ by \cite[Table~$5.4.C$]{KL}, and~$G$ is a
quasisimple group with an absolutely irreducible representation of 
degree~$R(T)$, we can use \cite{HM} to find candidates for~$G$ and thus~$S$. 
If $R( T ) = 248$, then $S = \text{\rm Th}$ is the only candidate, and this 
gives rise to case (f) (see \cite[p.~$177$]{Atlas} and \cite{KW}). Next, suppose 
that $R( T ) = 56$. Then $R( S ) < R( T )$ for all possible~$G$ by \cite{HM} and 
\cite{Jansen}. If $R( T ) = 27$, then $R( S ) < R( T )$ unless 
$S = \text{\rm Fi}_{22}$, $p = 2$ and $T = {^2\!E}_2( 2 )$. This gives rise to 
case (e) by \cite[Proposition~$5.4.17$]{KL} and \cite[p.~$162$]{Atlas}. If 
$R( T ) = 26$ or $R( T ) = 25$, then $R( S ) < R( T )$ for all possible~$G$.
If $R( T ) = 8$, then $T = {^3\!D}_4( p^a )$ for some positive integer~$a$. Here,
$R( S ) < R( T )$, unless $S = A_{10}$ and $p \in \{ 2, 5 \}$, 
or $S = \Omega^+_8( 2 )$ and $p$ is odd. Moreover, $X = \Sp_8( 2 )$ if $p = 2$, 
and $X = \Omega^+_8( p )$, if $p$ is odd. Neither of these groups~$X$ contains 
${^3\!D}_4( p^a )$ as a subgroup. Next, suppose that 
$R( T ) = 7$. Then $T = {^2G}_2(3^a)$, where~$a$ is an odd positive integer 
greater than~$1$. In particular, $p = 3$. This excludes the possibility $S = J_1$ 
by \cite{HM} and \cite{Jansen}. If $G = A_9$, then $X = \Omega_7( 3 )$, and 
$N_X( G )$ is maximal in~$X$ (see \cite[p.~$109$]{Atlas}). Finally, suppose that 
$R( T ) = 6$. Then $T = G_2( 2^a )$, where~$a$ is a positive integer greater 
than~$1$. In particular, $p = 2$. The tables in \cite{HM}
leave the options $S = J_2$ and $S = M_{22}$ (with $G = 3.M_{22}$). The first
case gives rise to case (c) by \cite[p.~$97$]{Atlas}, \cite{KW} and 
\cite[p.~$273$]{MAtl}. The second case yields $X = \SU_6( 2 )$ by 
\cite[p.~$93$]{MAtl} and the remark in \cite[p.~xi]{MAtl}, and $3.M_{22}$ is 
maximal in~$X$ by \cite[p.~$115$]{Atlas}.

Let~$\bar{G}$ denote the image of~$G$ in~$T$. Then~$\bar{G}$ is quasisimple 
and $\bar{G} \lneq T$.
Suppose now that $T$ is an alternating group of degree $m$. As $R( T ) \geq 5$,
we find $n = R( T ) = m - 2$ (see \cite[Proposition~$5.3.7$]{KL}), i.e., 
$T = A_{n + 2}$, and $p \mid n + 2$. The embedding $\bar{G} \rightarrow T$ 
yields a permutation representation of degree $n + 2$ of~$\bar{G}$. Let~$k$ be 
the size of a nontrivial orbit of~$\bar{G}$. Then~$\bar{G}$ has a subgroup of 
index~$k$, and thus some quasisimple quotient of~$\bar{G}$ embeds into $A_k$. 
By Remark~\ref{PropertiesR}(c) we find $R( A_k ) \geq R( S ) = n \geq 5$,
and hence $R( A_k ) = k - 2$. As $k \leq n + 2$, we conclude $k = n + 2$.
Also, the action of~$\bar{G}$ is two-fold transitive. Otherwise, the corresponding 
ordinary permutation character of~$\bar{G}$ would have a non-trivial constituent of 
degree smaller than~$n$, which would imply $R( S ) < n$, a contradiction. All 
possibilities 
for~$S$ and $n$ are listed in \cite[Table~$7.4$]{Cameron}. In each case 
$R( S ) < n$, as can be checked for sporadic groups $S$ from \cite{Jansen}.

Suppose finally that $T$ is a sporadic group. In this case we use the 
Atlas~\cite{Atlas} and~\cite[Section~$4$]{Wilson} to find the candidates 
for~$\bar{G}$ and hence~$S$, and \cite{Jansen} to rule out all possibilities 
except $S = \bar{G} = M_{23} \leq M_{24} = T$, $n = 11$ and $p = 2$. Using the 
remark in \cite[p.~xi]{MAtl} and the character tables in 
\cite[p.~$177$, p.~$267$]{MAtl}, we find that $X = \SL_{11}( 2 )$ and that 
$M_{23} \leq M_{24} \leq X$.
\end{prf}

\begin{rem}\label{SecondSmallest}
{\rm
If $G = X \lneq \SL( V )$ in Theorem~\ref{MaxEmb}, then $G$ is in Aschbacher 
class~$\mathcal{C}_8$. In this case, replace $X$ by the smallest classical 
group~$Y$ properly containing~$G$. In matrix notation, $Y = \Sp_n( q )$ if~$n$ 
and~$q$ are even and $G = \Omega^\pm_n( q )$. In all other cases 
$Y = \SL_n( q )$. Then $N_Y( G )$ is a maximal subgroup of~$Y$ (see 
\cite[Tables~$3.5$.A,~$3.5$.C]{KL} and \cite[Proposition~$2.3.32$]{BHRD}).
}
\end{rem}

\begin{rem}\label{SmallDegree}
{\rm
Assume Hypothesis~\ref{BasicHypothesis}, except for the condition that 
$n \geq 5$; suppose instead that $n \in \{2, 3, 4\}$. Then one of the
conclusions of Theorem~\ref{MaxEmb}~(a) or~(b) holds for~$G$. The statement of
Remark~\ref{SecondSmallest} is also true. 

Indeed, let~$S$ be a classical simple group of one of the following series:
$\PSL_n( q )$, $2 \leq n \leq 4$, $\PSU_n( q )$, $3 \leq n \leq 4$, 
$\PSp_n(q)$, $n = 4$, where $q$ is a power of the prime~$p$. Then 
$R( S ) = R_p( S ) = n$ and $R_\ell( S ) > n$ for all primes $\ell \neq p$,
unless 
\begin{equation}
\label{SmallGroups}
S \in \{ \PSL_2( 4 ) \cong \PSL_2( 5 ), \PSL_3( 2 ), \PSU_4( 2 ) \cong 
\PSp_4( 3 ) \}.
\end{equation}
Moreover, every finite simple group~$T$ with $R( T ) \leq 4$ is 
isomorphic to one of the classical groups listed above, to one of $A_m$, 
$5 \leq m \leq 8$, or to a Suzuki group. All these assertions follow from
\cite[Propositions~$5.4.13$,~$5.4.15$, Corollary~$5.4.14$]{KL} (if $T$ is a
sporadic group use~\cite{Jansen}). As the
alternating groups $A_m$, $5 \leq m \leq 8$, as well as the groups listed
in~(\ref{SmallGroups}) are contained in the Modular Atlas \cite{MAtl}, it
is easy to prove the initial statements using \cite[Tables~$8.1$--$8.17$]{BHRD}.
}
\end{rem}

\begin{rem}\label{MaximalOvergroups}
{\rm
Let $G \lneq K \lneq X$ be as in one of the cases (c)--(h) of Theorem~\ref{MaxEmb},
where~$K$ denotes the unique quasisimple group disproving the maximality of $N_X(G)$.
Then $N_X(K)$ is maximal in~$X$. Moreover, $N_X(K) = K$, except in cases (g) and (h),
where $N_X(K) = 3.{^2\!E}_6(2).3$ and $2 \times E_8( 3 )$, respectively.
Indeed, in cases (c) and (d), $\Aut( K )$ does not have an
irreducible representation of degree $R(K)$, as follows from
inspection of $\Aut( F_4(2) )$
and $\Aut( G_2( 3 ) )$, respectively. To prove the
statement in case (g), we use~\cite{WWWW} to obtain a $27$-dimensional 
irreducible representation of $3.{^2\!E}_6(2).3$ 
over $\mathbb{F}_4$, which can be shown with the help of the MeatAxe
\cite{MeatAxe} to preserve a hermitian form. In the remaining cases,
$\Aut(K) = K$ so that $N_X(K) = C_X( K ) = KZ(X)$.
}
\end{rem}

\section*{Acknowledgements} 
I thank Frank L{\"u}beck for helpful conversations and Thomas Breuer for his 
careful reading of a first version of this manuscript. I also gratefully 
acknowledge support by the German Research Foundation (DFG) within the 
SFB-TRR 195 ``Symbolic Tools in Mathematics and their Application'', to which 
this work is a contribution. Finally, I am very much indebted to the referee 
of a first version of this paper for pointing out the omission of the Ree 
groups in the statement of the main theorem.

\end{document}